\documentclass[11pt]{amsproc}
\usepackage{
amssymb,
amsmath}

\usepackage[bookmarks,colorlinks,pagebackref]{hyperref} 

\usepackage{mathptmx}
\usepackage{helvet}
\usepackage{courier}

\usepackage[all,cmtip,2cell]{xy}

\numberwithin{equation}{section}
\theoremstyle{plain} 
\newtheorem{proposition}{Proposition}[section]  
\newtheorem{lemma}[proposition]{Lemma}
\newtheorem{corollary}[proposition]{Corollary} 
\newtheorem{theorem}[proposition]{Theorem} 

\theoremstyle{definition} 
\newtheorem{definition}[proposition]{Definition}

\newtheorem{remark}[proposition]{Remark} 

\newtheorem{example}[proposition]{Example} 

\newtheorem{notation}[proposition]{Notation}

\newtheorem{notations}[proposition]{Notations}

\newcommand\Supp{\operatorname{Supp}}

\newcommand\Ann{\operatorname{Ann}}

\newcommand\Tor{\operatorname{Tor}}
\newcommand\Hom{\operatorname{Hom}}

\newcommand\Ext{\operatorname{Ext}}
\newcommand\Rad{\operatorname{Rad}}
\newcommand\Ker{\operatorname{Ker}}
\newcommand\Coker{\operatorname{Coker}}

\newcommand{\xx}{\underline x}
\newcommand{\yy}{\underline y}

\newcommand\grade{\operatorname{grade}}

\newcommand\cd{\operatorname{cd}}

\newcommand\depth{\operatorname{depth}}
\newcommand\height{\operatorname{height}}
\newcommand\Spec{\operatorname{Spec}}
\newcommand{\gam}{\Gamma_{I}}
\newcommand{\Char}{\operatorname{char}}
\newcommand{\Rgam}{{\rm R} \Gamma_I}

\newcommand{\Llam}{{\rm L}\Lambda^I}

\newcommand{\qism}{\stackrel{\sim}{\longrightarrow}}

\author[P.~Schenzel]{Peter Schenzel}
\title[Endomorphisms]
{Notes on endomorphisms, local cohomology and completion}
\address{Martin-Luther-Universit\"at Halle-Wittenberg,
Institut f\"ur Informatik, D --- 06 099 Halle (Saale), Germany}
\email{peter.schenzel@informatik.uni-halle.de}

\begin{document}

\begin{abstract} 
Let $M$ denote a finitely generated module over a Noetherian ring $R$. 
For an ideal $I \subset R$ there is a study of the endomorphisms of the 
local cohomology module $H^g_I(M), g = \grade (I,M),$ and related results. 
Another subject 
is the study of left derived functors of the $I$-adic completion $\Lambda^I_i(H^g_I(M))$, motivated by 
a characterization of Gorenstein rings given in \cite{SS}.	This provides another 
Cohen-Macaulay criterion. The results are illustrated by several examples. There is also an extension to the case of homomorphisms of two different local cohomology modules.
\end{abstract}


\maketitle
\begin{center}
	\textsl{In honour of Sylvia and Roger Wiegand for their 
	contributions in Commutative Algebra.}
\end{center}


\section{Introduction}
Let $I$ denote an ideal of a Noetherian ring $R$. Let $M$ 
be a finitely generated $R$-module, and let $H^i_I(M), i \in \mathbb{Z},$ 
denote the local cohomology modules of $M$ with respect to $I$ (see 
\cite{Ga2} or \cite{BrS} for definitions). In the case of $(R,\mathfrak{m})$, a local ring, we denote by $E_R(\Bbbk)$ 
the injective hull of the residue field $R/\mathfrak{m} = \Bbbk$ and by 
$\Hom_R(\cdot,E_R(\Bbbk)) = D(\cdot)$ the Matlis Duality functor. 

In their paper (see \cite{HS}) 
Hellus and St\"uckrad investigated the endomorphism ring of $H^g_I(R)$ 
and $D(H^g_I(R))$ in the case of a complete local ring $R$ and an ideal 
satisfying $H^i_I(R) = 0$ for all $i \not= g$. In fact, they proved the isomorphisms of endomorphism rings 
\[
R \cong \Hom_R(H^g_I(R),H^g_I(R)) \cong \Hom_R(D(H^g_I(R)),D(H^g_I(R))).
\] 
A characterization of ideals such that $H^i_I(R) = 0$ for all $i \not= g$, 
so-called cohomologically complete intersections, is described by Hellus 
and the author (see \cite{HSp2}). Some of these results were generalized 
to finitely generated modules by Zargar (see \cite{Zmr}) under the name 
relative Cohen-Macaulay modules. In the case of a Gorenstein ring $(R,\mathfrak{m})$ 
weaker conditions than $H^i_I(R) =0$ for $i \not= 
\grade (I,R)$ are sufficient for the previous isomorphisms (see \cite{Sp6}).

In the following we shall extend some of the previous results to the 
situation of an $R$-module $M$ and extend the result by Hellus and St\"uckrad (see \ref{cor-3}). Namely we will investigate the local cohomology 
module $H^g_I(M), g = \grade (I,M)$. To this end we denote by $\hat{(\cdot)}{}^I$ the 
$I$-adic completion functor. For the maximal ideal $I = \mathfrak{m}$ we write 
$\hat{(\cdot)}$. As main results we shall prove the following:

\begin{theorem} \label{thm-4}
	Let $I$ denote an ideal of a local ring $(R, \mathfrak{m})$. Let $M$ be a finitely 
	generated $R$-module and $g = \grade (I,M)$. 
	\begin{itemize}
		\item[(a)] $\Hom_{\hat{R}}(H^g_{I\hat{R}}(\hat{M}),H^g_{I\hat{R}}(\hat{M})) \cong \Hom_R(D(H^g_I(M)), D(H^g_I(M)))$ and therefore
		$$
		\hat R \cong \Hom_{\hat R}(H^g_{I \hat R}(\hat{M}),H^g_{I \hat R}(\hat{M})) \mbox{ if and only if } \hat R \cong \Hom_R(D(H^g_I(M)), D(H^g_I(M))).
		$$
		\item[(b)] If $H^i_I(M) = 0$ for all $i \not= g$, then 
		$$
		\Hom_R(M,M) \otimes_R \hat{R}^I 
		\cong \Hom_R(H^g_I(M),H^g_I(M)) \mbox{ and }	\Ext_R^j(M,M) \otimes_R \hat{R}^I \cong \Ext_R^{g+j}(H^g_I(M),M) \mbox{ for all } j.
		$$
	\end{itemize}
\end{theorem}

Extensions of the isomorphisms in \ref{thm-4} (b) for two different modules 
is given in Section 6.
Moreover, there is also another characterization of the grade of a module. 
For our further results we denote by $\Lambda_i^I(\cdot)$ the left derived functor 
of the completion $\hat{(\cdot)}{}^I$. We refer to \cite{SS} for definitions and 
basic results. In \cite[10.5.9]{SS} it is shown that $R$ is a $d$-dimensional 
Gorenstein ring if and only if $\Lambda_i^{\mathfrak{m}}(E_R(\Bbbk)) = 0$ for 
all $i \not= d$ and $\Lambda_d^{\mathfrak{m}}(E_R(\Bbbk)) \cong \hat{R}$. As a 
generalization we prove the following result:

\begin{theorem} \label{thm-5}
	Let $I \subset R$ be an ideal. Let $M$ denote a finitely generated $R$-module 
	and $g = \grade (I,M)$.
	\begin{itemize}
		\item[(a)] There is a natural homomorphism $\Lambda_g^I(H^g_I(M))  \stackrel{\tau}{\longrightarrow} \hat{M}^I$.
		\item[(b)] Suppose  $H^i_I(M) = 0$ for all $i \not= g$, then $\tau$ 
		is an isomorphism and $\Lambda_i^I(H^g_I(M)) = 0$ for all $i \not= g$.
	\end{itemize}
\end{theorem}

As an application it yields a further Cohen-Macaulay criterion: $M$ is a 
$d$-dimensional Cohen-Macaulay module if and only if 
$\Lambda_d^{\mathfrak{m}}(H^d_{\mathfrak{m}}(M)) \cong \hat{M}^{\mathfrak{m}}$ 
and $\Lambda_i^{\mathfrak{m}}(H^d_{\mathfrak{m}}(M)) = 0$ for all $i \not=d$ 
(see \ref{cor-4}). Furthermore, there is an intrinsic characterization 
of canonically Cohen-Macaulay modules in the sense of \cite{Sp4}. A few more results about modules satisfying the assumption in \ref{thm-4} (b) are shown.

In our terminology we follow the textbooks \cite{BrS} and \cite{Mh}. As a reference 
for Homological Algebra we refer to \cite{EJ}. Particular results about 
local duality and completion as well as its derived functors are found in \cite{SS}. By "$\simeq$" we denote quasi-isomorphisms in the sense of \cite[1.1.3]{SS}.
In a final section we illustrate the results by several examples. In particular, 
$\Hom_R(H^g_I(M),H^g_I(M), g = \grade(I,M),$ is not a finitely generated $R$-module in general (see \ref{expl-4}).

\section{Preliminaries and the truncation complex}
\begin{notation} \label{not-1}
(A) In the following let $R$ denote a commutative $d$-dimensional 
ring. If $(R,\mathfrak{m})$ is local, let $\Bbbk = R/\mathfrak{m}$ be its residue field. In this case we denote by 
$E = E_R(\Bbbk)$ the injective hull of the residue field and by $D(\cdot) = 
\Hom_R(\cdot, E)$ the Matlis Duality functor. Let $I \subset R$ be an ideal 
and $M$ be an $R$-module. For a prime ideal $\mathfrak{p} \in \Spec R$ let 
$E_R(R/\mathfrak{p})$ denote the injective hull of $R/\mathfrak{p}$ as $R$-module.\\
(B) Let $M$ denote an $R$-module. Let $E_R^{\cdot}(M) : 0 \to M \to E_R^0(M) \to E_R^1(M) \to \ldots \to 
E_R^i(M) \to \ldots$ be a minimal injective resolution. By Matlis' Structure 
theory it follows that $E_R^i(M) \cong \oplus_{\mathfrak{p} \in \Supp_R M} 
E_R(R/\mathfrak{p})^{\mu_i(\mathfrak{p},M)}$, where $\mu_i(\mathfrak{p},M) = 
\dim_{k(\mathfrak{p})} \Ext_{R_{\mathfrak{p}}}(k(\mathfrak{p}),M_{\mathfrak{p}})$.
For the details we refer to \cite{Bh} and \cite{EJ}. \\
(C) For a finitely generated $R$-module $M$ we define (following D. Rees) 
$$
\grade(I,M) = \min \{ i \in \mathbb{N}| \Ext_R^i(R/I,M) \not= 0\}.
$$ 
For basic properties of $\grade(I,M)$ we refer to \cite[1.2]{BrH}, in particular 
$\grade(I,M) = \inf \{ \depth M_{\mathfrak{p}}| \mathfrak{p} \in V(I) \}$ 
(see e.g. \cite[1.2.10]{BrH}).\\
(D) With the previous notation let $H^i_I(M), i \in \mathbb{N},$ denote the 
local cohomology modules of $M$ with respect to $I$. We refer to \cite{Ga2} 
or to \cite{SS} for some basic results. Note that $H^i_I(\cdot)$ is the $i$-th 
derived functor of the section functor $\Gamma_I(\cdot)$. 
\end{notation}

In the following we shall modify the construction of the truncation complex 
as introduced by the author (see \cite{Sp1} and \cite{Sp4}) and used several times. 

\begin{definition} \label{def-1} {\textsl{Truncation complex.}}
	With the previous notation let $g = \grade(I,M)$. We apply the section 
	functor $\Gamma_I(\cdot)$ to the minimal injective resolution $E_R^{\cdot}(M)$. 
	By view of \ref{not-1} (B) it follows that $\Gamma_I(E_R^{\cdot}(M))^i = 0$ 
	for all $i < g$. So there is 
	an injection $0 \to H^g_I(M)[-g] \to \Gamma_I(E_R^{\cdot}(M))$, where 
	$H^g_I(M)[-g]$ denotes the module $H^g_I(M)$ considered as a complex sitting 
	in cohomological degree $g$. The complex $C_M^{\cdot}(I)$ obtained as the 
	cokernel is denoted  the truncation complex. So there is a short exact sequence of complexes 
	\[
	0 \to H^g_I(M)[-g] \to \Gamma_I(E_R^{\cdot}(M)) \to C_M^{\cdot}(I) \to 0
	\]
	and $H^i(C_M^{\cdot}(I)) \cong H^i_I(M)$ for $i > g$ and $H^i(C_M^{\cdot}(I)) 
	= 0$ for $i \leq g$.
\end{definition}

Next we shall prove some homological results related to modules with the 
property that their support is contained in $V(I)$. 

\begin{lemma} \label{lem-1} Let $I \subset R$ denote an ideal. 
	Let $M, N$ denote two $R$-module. Suppose that $\Supp_R M \subseteq V(I).$ 
	There are the following natural isomorphisms:
	\begin{itemize}
		\item[(a)] $\Hom_R(M, \Gamma_I(N)) \simeq \Hom_R(M,N),$
		\item[(b)] $M \otimes_R \Hom_R(\Gamma_I(N),E_R(\Bbbk)) \simeq M \otimes_R \Hom_R(N,E_R(\Bbbk)).$
	\end{itemize}
\end{lemma}
\begin{proof}
	The proof in (a) is easy to see. For the proof of (b) assume at first that $M$ 
	is a finitely generated $R$-module. Then $\Hom_R(\Hom_R(M,\cdot),E_R(\Bbbk) 
	\cong M\otimes_R \Hom_R(\cdot,E_R(\Bbbk))$ and the claim follows by 
	(a) and adjointness. The general case turns out because any $R$-module 
	can be written as diret limit of finitely generated submodules and tensor products commute 
	with direct limits.
\end{proof}

\begin{lemma} \label{lem-2}
	Let $X, M$ be  $R$-modules and $\Supp_R X \subseteq V(I),$
	where $I \subset R$ denotes an ideal with $g = \grade (I,M).$ 
	\begin{itemize}
		\item[(a)] $\Ext^i_R(X,M) = 0$ for $i < g$ and $\Ext^g_R(X,M) \cong \Hom_R(X, H^g_I(M)).$
		\item[(b)] $\Tor^R_i(X, D(M)) = 0$ for $i < g$ and $\Tor_g^R(X, D(M)) \cong X \otimes_R D(H^g_I(M)).$
	\end{itemize}
\end{lemma}

\begin{proof}
	As above let $M \qism E_R^{\cdot}(M)$ denote a minimal injective resolution.
	By view of Lemma \ref{lem-1} there is an isomorphism of
	complexes
	\[
	\Hom_R(X, \Gamma_I(E^{\cdot}_R(M))) \cong \Hom_R(X, E^{\cdot}_R(M))
	\]
	since $\Supp X \subseteq V(I).$
	By \ref{def-1} there is an exact sequence $0 \to H^g_I(M) \to \Gamma(E^{\cdot}(M))^g \to \Gamma(E^{\cdot}M))^{g+1}.$
	Because of the previous remark it induces a natural commutative diagram with exact rows
	\[
	\begin{array}{cccccc}
	0 \to & \Hom_R(X, H^g_I(M)) & \to & \Hom_R(X,\gam(E^{\cdot}(M)))^g  & \to  & \Hom_R(X,\gam(E^{\cdot}(M)))^{g+1} \\
	&     \downarrow             &      & \downarrow                          &      & \downarrow \\
	0 \to & \Ext_R^g(X, M)        & \to  & \Hom_R(X, E^{\cdot}(M))^g       & \to  & \Hom_R(X, E^{\cdot}(M))^{g+1}.
	\end{array}
	\]
	The last two vertical homomorphisms are isomorphisms by Lemma \ref{lem-1}. 
	Therefore the first vertical map
	is an isomorphism too. This finishes the proof of (a). 
	
	For the proof of (b) 
	we first assume that $X$ is a finitely generated $R$-module. By the 
	$\Ext-\Tor$ duality (see e.g. \cite[1.4.1]{SS}) there are canonical isomorphisms
	\[
	\Hom_R(\Ext^i_R(X,M),E_R(\Bbbk)) \cong \Tor_i^R(X,\Hom_R(M,E_R(\Bbbk)))
	\]
	for all $i$. 
	This proves the statement (b) for a finitely generated $R$-module $X$. 
	In general let $X = \varinjlim X_{\alpha}$ be the direct limit of $X$ by finitely 
	generated submodules $X_{\alpha} \subset X$. Because the tensor product  
	and $\Tor$ commute with direct limits this completes the proof.
\end{proof}

\begin{remark} \label{rem-1}
	The previous result \ref{lem-2} (b) provides another characterization 
	of $\grade(I,M)$ for a finitely generated module $M$ over a local ring. Namely 
	$\grade(I,M) = \min \{i \in \mathbb{N} | \Tor_i^R(R/I,D(M)) \not= 0 \}$. 
	This is true because of  $\Tor_i^R(R/I,D(M)) \cong D(\Ext_R^i(R/I,M))$ for all $i$. 
	The dual notion of the grade is the $\operatorname{witdth} (I,M)$ as introduced by 
	Frankild (see \cite{Fa}). It is related to the minimal non-vanishing of local homology. See also \cite[5.1.2]{SS}, where it is denoted by Tor-codepth. 
\end{remark}

\section{Endomorphisms of Local cohomology}
The statements of Lemma \ref{lem-2} provide some information about 
the endomorphisms of the local cohomology module $H^g_I(M), g = \grade(I,M)$. 
In the following we will denote by $\hat{(\cdot)}{}^I$ the $I$-adic completion 
functor. In the case of a local ring $(R,\mathfrak{m})$ we write $\hat{(\cdot)}$ 
for the $\mathfrak{m}$-adic completion functor.

\begin{theorem} \label{thm-1}
	Let $I$ denote an ideal of a local ring $(R,\mathfrak{m})$. Let $M$ be a 
	finitely generated $R$-module with $g = \grade(I,M)$. Then there is a natural 
	isomorphism 
	\[
	\Hom_{\hat{R}}(H^g_{I\hat{R}}(\hat{M}),H^g_{I\hat{R}}(\hat{M})) \cong 
	\Hom_R(D(H^g_I(M)),D(H^g_I(M))).
	\]
\end{theorem}

\begin{proof}
	By adjunction formula we get the natural isomorphism 
	\[
	\Hom_R(D(H^g_I(M)),D(H^g_I(M))) \cong D(H^g_I(M) \otimes_R D(H^g_I(M))).
	\]
	By virtue of \ref{lem-2} (b) there are the following isomorphisms to the second 
	module above
	\[
	D(\Tor_g^R(H^g_I(M),D(M) ) \cong 
	\Ext_R^g(H^g_I(M), \hat{M}) \cong \Ext^g_{\hat{R}}(H^g_{I\hat{R}}(\hat{M}),\hat{M}).
	\]
	For the first recall the $\Ext-\Tor$ duality (see e.g. \cite[1.4.1]{SS}) 
	and $D(D(M)) \cong \hat{M}$. For the second we refer to 
	\cite[Lemma 2.1]{Sp5} and recall that 
	$H^g_I(M) \otimes_R \hat{R} \cong H^g_{I\hat{R}}(\hat{M})$. Then the 
	statement follows by \ref{lem-2} (a) applied to $\hat{R}$ and $\hat{M}$.
\end{proof}

For an $R$-module $X$ there is a natural $R$-homomorphism $\phi: R \to \Hom_R(X,X), 
r \mapsto rx$, the multiplication map by $x$ on $X$. Clearly $\ker \phi = \Ann_R X$.
In the following we shall have a look for the situation of $X = H^g_I(M)$. 
To this end we look also at the natural map $X \otimes_R \Hom_R(X,Y), 
x \otimes \psi \mapsto  \psi(x)$, where $X,Y$ denote $R$-modules. 

\begin{corollary} \label{cor-1}
	With the previous notation  the following conditions about the corresponding 
	natural $R$-homomorphisms are equivalent:
	\begin{itemize}
		\item[(i)] $H^g_I(M) \otimes D(H^g_I(M)) \to E_R(\Bbbk)$ is an isomorphism.
		\item[(ii)] $\hat R \to \Hom_{\hat R}(H^g_{I \hat R}(\hat{M}),H^g_{I \hat R}(\hat{M}))$ is an isomorphism.
		\item[(iii)] $\hat R \to \Hom_R(D(H^g_I(M)), D(H^g_I(M)))$ is an isomorphism.
	\end{itemize}
\end{corollary}

\begin{proof}
	The equivalence of (ii) and (iii) follows by Theorem \ref{thm-1}. 
	The natural homomorphism in (i) is an isomorphism if and only if $D(E_R(\Bbbk)) \to 
	D(H^g_I(M) \otimes D(H^g_I(M)))$ is an isomorphism. The last module 
	is isomorphic to $\Hom_R(D(H^g_I(M)), D(H^g_I(M)))$. It is easy to see that this 
	is the natural map.
\end{proof}

Note that \ref{thm-1} and \ref{cor-1} prove Theorem \ref{thm-4} (a). The 
following comparison result could be of some interest in order to relate 
different endomorphism rings to each other. 

\begin{remark} \label{rem-2}
Let $J \subseteq I$ denote 
two ideals such that $\grade(I,M) = \grade(J,M) =g$. Then there is an exact sequence 
\[
0 \to H^g_I(M) \stackrel{\psi}{\to} H^g_J(M) \to 
H^{g+1}_{V(J)\setminus V(I)}(\tilde{M}) 
\]
as easily seen. Here $\tilde{M}$ denotes the associated sheaf of $M$ on $\Spec R$. 
Let $X = \Coker \psi$ then it induces an exact sequence 
\[
0 \to \Hom_R(H^g_J(M),H^g_J(M)) \to \Hom_R(H^g_I(M),H^g_I(M)) \to \Ext_R^{g+1}(X,M)
\]
as follows by applying $\Hom_R(\cdot,M)$ and by view of \ref{lem-2}. Note that 
$\Supp_RX \subseteq V(J)$ and use \ref{lem-2} (a). 
\end{remark}

In the following we shall investigate when $\psi$ is an isomorphism. 

\begin{lemma} \label{lem-3}
	Let $J \subseteq I$ denote two ideals with $\grade(I,M) = \grade(J,M) =g$. 
	Suppose that $\depth_{R_{\mathfrak p}}M_{\mathfrak p} \geq g+1$ for all 
	$\mathfrak p \in \Supp_R M \cap (V(J) \setminus V(I))$.  Then the homomorphisms
		\[
		H^g_I(M)\to H^g_J(M) \mbox{ and }
		\Hom_R(H^g_J(M), H^g_J(M)) \to \Hom_R(H^g_I(M), H^g_I(M))
		\]
		are isomorphisms.
\end{lemma}

\begin{proof}
	For the proof it will be enough to prove the first isomorphism. 
	To this end recall that $\Supp_RX \subseteq \Supp_R M \cap (V(J) \setminus V(I))$. 
	By the assumption $\depth_{R_{\mathfrak p}} M_{\mathfrak p} \geq g+1$ for all $\mathfrak{p} 
	\in \Supp_R M \cap (V(J) \setminus V(I))$. Whence $H^{g+1}_{V(J)\setminus V(I)}(\tilde{M}) = 0$ 
	and $\psi$ is an isomorphisms. 
\end{proof}

Note that $\depth_{R_{\mathfrak p}}M_{\mathfrak p} \geq g$ for all 
$\mathfrak p \in \Supp_R M \cap (V(J) \setminus V(I))$ as follows by the 
characterization of $\grade(J,M) = g$. Moreover $\Supp_R M \cap (V(J) \setminus V(I)) = 
\Supp_R M/JM :_M\langle I \rangle$, where $JM :_M\langle I \rangle 
$ denotes the stable value of $JM :_M I^n$ for $n \gg 0$.

\section{Endomorphisms of modules}

In this section we shall relate the endomorphism ring of an $R$-module $M$ to 
that of a certain local cohomology module. An approach to the 
$\Ext$-cohomology of $H^g_I(M)$ is given in the following result.

\begin{theorem} \label{thm-2}
	With the above notation let $M, N$ be a finitely generated $R$-modules  
	and $g = \grade(I,M)$.
	Let $E_R^{\cdot}(M), E_R^{\cdot}(N)$ be  minimal injective resolutions of $M$
	and $N$ resp. Then there is a 
	long exact sequence 
	\begin{gather*}
	\ldots \to H^i(\Hom_R(C^{\cdot}_M(I),E_R^{\cdot}(N))) \to 
	\Ext_R^i(M,N) \otimes_R \hat{R}^I \to \Ext_R^{g+i}(H^g_I(M),N) \\ 
	\to H^{i+1}(\Hom_R(C^{\cdot}_M(I),E_R^{\cdot}(N))) \to \ldots 
	\end{gather*}
\end{theorem}

\begin{proof}
	We start with the short exact sequence of the truncation complex as given 
	in \ref{def-1}. Then we apply the functor $\Hom_R(\cdot,E_R^{\cdot}(N))$. 
	Therefore there is the following short exact sequence of complexes
	\[
	0 \to \Hom_R(C^{\cdot}_M(I),E_R^{\cdot}(N)) \to \Hom_R(\Gamma_I(E_R^{\cdot}(M)),E_R^{\cdot}(N)) \to 
	\Hom_R(H^g_I(M),E_R^{\cdot}(N)[g] \to 0
	\]
	because $E_R^{\cdot}(N)$ is a left bounded complex of injective modules.
	The $i$-th cohomology of the last complex is $\Ext_R^{g+i}(H^g_I(M),N)$. 
	Next we inspect the cohomology of the complex in the middle. To this end 
	let $\xx = x_1,\ldots,x_t$ denote a generating set of the ideal $I$ and let 
	$\check{L}_{\xx}$ denote the bounded free resolution of the $\check{C}_{\xx}$ 
	complex as constructed in \cite[6.2.2 and 6.2.3]{SS} or \cite{Sp2}. 
	Note that $M \to E_R^{\cdot}(M)$ is a quasi-isomorphism 
	and therefore $\Gamma_I (E_R^{\cdot}(M)) \simeq \check{L}_{\xx} \otimes_R M$. 
	Because 
	$E_R^{\cdot}(N)$ is a left bounded complex of injective modules the complex 
	in the middle is quasi-isomorphic to 
	\[
	\Hom_R(\check{L}_{\xx} \otimes_R M, E_R^{\cdot}(N)) \cong \Hom_R(\check{L}_{\xx}, 
	\Hom_R(M,E_R^{\cdot}(N))),
	\] 
	(see also \cite{SS} for more details).  Now it follows that 
	\begin{equation}
	H_i(\Hom_R(\check{L}_{\xx}, \Hom_R(M,E_R^{\cdot}(N)))) 
	\cong \Lambda_i^I(\Hom_R(M,E_R^{\cdot}(N))). \tag{$\star$}
	\end{equation}
	Since the cohomology $\Ext_R^i(M,N)$, the cohomology of $\Hom_R(M,E_R^{\cdot}(N)$, 
	is finitely generated it implies the isomorphisms 
	\[
	H^i(\Hom_R(\Gamma_I(E_R^{\cdot}(M)),E_R^{\cdot}(N))) \cong \Ext_R^i(M,N) \otimes_R \hat{R}^I
	\] 
	for all $i$. By the long exact cohomology sequence it proves the claim.
\end{proof}

\begin{corollary} \label{cor-2}
	With the above notation we get the following:
	\begin{itemize}
		\item[(a)] There is a natural homomorphism $\Hom_R(M,M) \otimes_R \hat{R}^I 
		\to \Hom_R(H^g_I(M),H^g_I(M)))$.
		\item[(b)] Suppose that $H^i_I(M) = 0$ for all $i \not= g$. Then the 
		map in (a) is an isomorphism and 
		\[
		\Ext_R^j(M,M) \otimes_R \hat{R}^I \cong \Ext_R^{g+j}(H^g_I(M),M)
		\]
		for all $j$.
	\end{itemize}
\end{corollary}

\begin{proof}
	By \ref{lem-2} (b) there is an isomorphism $\Hom_R(H^g_I(M),H^g_I(M))) 
	\cong \Ext_R^g(H^g_I(M),M)$. Whence the statement in (a) follows \ref{thm-2}. 
	
	Under the additional assumption in (b)
	$C^{\cdot}_M(I)$ is homologically trivial 
	and $\Hom_R(C^{\cdot}_M(I),E_R^{\cdot}(M))$ is exact. Therefore the 
	statements follow by (a) and \ref{thm-2}.
\end{proof}

Note that \ref{cor-2} proves \ref{thm-4} (b). As a further application we 
get a result concerning a cohomological complete intersection, i.e. an ideal 
such that $H^i_I(R) = 0$ for all $i \not= \grade(I,R)$. 

\begin{proposition} \label{prop-1}
	Let $I \subset R$ an ideal with $H^i_I(R) = 0$ for all $i \not= \grade(I,R)$.
	Let $N$ denote an arbitrary $R$-module. Then there are isomorphisms
	\[
	\Ext_R^{g-i}(H^g_I(R),N) \cong \Lambda_i^I(N) 
	\] 
	for all $i \in \mathbb{Z}$.
\end{proposition}

\begin{proof}
	We follow a slight modification of the proof of \ref{thm-2}. By the short exact sequence at the beginning of the proof of \ref{thm-2} it follows that 
	$\Hom_R(\Gamma_I(E_R^{\cdot}(R)),E_R^{\cdot}(N)) \qism 
	\Hom_R(H^g_I(R),E_R^{\cdot}(N)[g]$ since $C_R^{\cdot}(I)$ is homologically trivial. For the 
	isomorphisms $(\star)$ in the proof of \ref{thm-2} with $M = R$ it yields 
	\[
	H_i(\Hom_R(\check{L}_{\xx},E_R^{\cdot}(N))) 
	\cong \Lambda_i^I(N).
	\]
	For its proof in \ref{thm-2} it is not necessary to assume that $N$ is 
	finitely generated. Because $E_R^{\cdot}(N)$ is an injective resolution 
	we have $\Hom_R(\check{L}_{\xx},E_R^{\cdot}(N))) \simeq \Hom_R(\check{L}_{\xx},N)$ and therefore $H_i(\Hom_R(\check{L}_{\xx},N)) 
	\cong \Lambda_i^I(N)$. Then the claim follows by the previous quasi-isomorphism 
	(see also \ref{not-2}).
\end{proof}

Another application is a result related to \ref{cor-1}. This is an
improvement of the result of Hellus and St\"uckrad (see \cite{HS}). For the proof we refer to \ref{cor-2} and \ref{prop-1}. 
For further results on the endomorphism ring of $H^g_I(R), g = \grade (I,R),$ 
under additional assumptions on the ring we refer to \cite{Sp6}. 

\begin{corollary} \label{cor-3}
	Let $I$ denote an ideal of a Noetherian ring $R$ and 
	$g = \grade(I,R)$. Suppose that $H^i_I(R) = 0$ for all 
	$i \not= g$. 
	\begin{itemize}
		\item[(a)] There are isomorphisms $\Ext_R^i(H^g_I(R),H^g_I(R)) 
		\cong \Lambda_{g-i}^I(H^g_I(R))$ for all $i \in \mathbb{Z}$. 
		\item[(b)] The canonical map $\hat{R}^I \to \Hom_R(H^g_I(R),H^g_I(R)))$
		is an isomorphism, $\Hom_R(H^g_I(R),H^g_I(R)) \cong \Lambda_g^I(H^g_I(R))$ 
		and $\Ext_R^i(H^g_I(R),R) = 0$ for all $i \not= g$.
	\end{itemize}
\end{corollary}

The previous results have an interesting application on the annihilators 
of certain local cohomology modules. 

\begin{remark} \label{rem-3}
	Let $I$ denote an ideal in a complete local ring  $(R,\mathfrak{m})$. Let $M$ 
	be  a finitely generated $R$-module $M$ such that $H^i_I(M) = 0$ for all 
	$i \not= \grade(I,M)$. Then $\Ann_RM = \Ann_R H^g_I(M)$. This follows since the natural 
	homomorphism $\Phi: R \to \Hom_R(X,X)$ as in \ref{cor-1} is given by 
	$r \mapsto \phi(r)$ with 
	$\phi(r) :X \to X, x \mapsto rx,$ the multiplication map and since $\Ker \Phi = \Ann_R X$.
\end{remark}

\section{Relation to completions}
\begin{notations} \label{not-2}
	(A) Let $\xx = x_1,\ldots, x_t$ denote a system of elements in $(R,\mathfrak{m})$ 
	and $I = (\xx)R$. The we recall the definition of the \v{C}ech $\check{C}_{\xx}$ 
	and its free resolution $\check{L}_{\xx}$ as done in 
	\cite[6.2.2 and 6.2.3]{SS} resp. \cite{Sp2}.\\
	(B) Let $M$ be an $R$-module. In the case of a Noetherian ring it 
	follows that $\Lambda_i^I(M) \cong H_i(\Hom_R(\check{L}_{\xx},M))$ 
	(see \cite{SS} for more details and generalizations). \\
	(C) Here $\Lambda_i^I(M), i \geq 0,$ denote the $i$-th left derived functor 
	of the $I$-adic completion with respect to $M$.
\end{notations}

In the following we shall provide another application of the truncation 
complex.

\begin{theorem} \label{thm-3}
	Let $I \subset R$ denote an ideal. 
	Suppose that $M$ is a finitely generated $R$-module and $g = \grade(I,M)$.
	Then there are an  exact sequence and isomorphisms
	\begin{gather*}
	0\to \Lambda_1^I(C_M^{\cdot}(M)) \to \Lambda_g^I(H^g_I(M))  \stackrel{\tau}{\longrightarrow} \hat{M}^I 
	\to \Lambda_0^I(C_M^{\cdot}(M)) \to \Lambda_{g-1}^I(H^g_I(M)) \to 0, \\
	\Lambda_{g+i}^I(H^g_I(M)) \cong \Lambda_{i+1}^I(C_M^{\cdot}(M)) \mbox{ for all } i \geq 1 \mbox{ and } i < -1.	
	\end{gather*}
	Suppose that $H^i_I(M) = 0$ for all $i \not= g$ then $\tau$ is an isomorphism 
	and $\Lambda_j^I(H^g_I(M)) = 0$ for all $j \not= g$.
\end{theorem}

\begin{proof}
	Let $E_R^{\cdot}(M)$ denote a minimal injective resolution of $M$. Then we use 
	the truncation complex (see \ref{not-1}). We use the notions of the proof of 
	\ref{thm-1} and we apply $\Hom_R(\check{L}_{\xx},\cdot)$ and obtain the following 
	short exact sequence 
	of complexes
	\[
	0 \to \Hom_R(\check{L}_{\xx},H^g_I(M))[-g] \to \Hom_R(\check{L}_{\xx},\Gamma_I(E_R^{\cdot}(M))) 
	\to \Hom_R(\check{L}_{\xx},C_M^{\cdot}(I)) \to 0.
	\]
	Recall that $\check{L}_{\xx}$ is a bounded complex of free $R$-modules. 
	Next we investigate the complex in the middle. By \cite[Chapter 7]{SS} 
	note the following quasi-isomorphisms 
	\[
	\Gamma_I(E_R^{\cdot}(M)) \simeq \check{L}_{\xx} \otimes_R E_R^{\cdot}(M) 
	\simeq \check{L}_{\xx} \otimes_R M.
	\]
	Therefore the complex in the middle is quasi-isomorphic to 
	$\Hom_R(\check{L}_{\xx},\check{L}_{\xx}\otimes_RM) \simeq \Hom_R(\check{L}_{\xx},M)$, 
	where the last quasi-isomorphism follows by view of \cite[6.5.4]{SS}. Now 
	\[
	H_i(\Hom_R(\check{L}_{\xx},M)) \cong \Lambda_i^I(M) \mbox{ for all } i \geq 0
	\] 
	(see \ref{not-2} (B)). Because $M$ is finitely generated it follows that 
	$\Lambda_0^I(M) \cong \hat{M}^I$ and the vanishing for $i \not= 0$. 
	
	That is, the long exact homology sequence provides the first part of the 
	statement. If $H^i_I(M) = 0$ for all $i \not= g$ then $C_M^{\cdot}(I)$ and 
	also $\Hom_R(\check{L}_{\xx},C_M^{\cdot}(I))$ is exact. Therefore,  
	by the long exact cohomology sequence it finishes the proof.
\end{proof}

The previous Theorem proves \ref{thm-4} of the Introduction. In the following 
corollary we will discuss the situation of a Cohen-Macaulay module. 

\begin{corollary} \label{cor-4}
	Let $M$ denote a $d$-dimensional  module over a local ring 
	$(R,\mathfrak{m})$. Then the following conditions are equivalent:
	\begin{itemize}
		\item[(i)] $M$ is a Cohen-Macaulay module.
		\item[(ii)] $\Lambda_d^{\mathfrak{m}}(H^d_{\mathfrak{m}}(M)) 
		\cong \hat{M}^{\mathfrak{m}}$ 
		and $\Lambda_i^{\mathfrak{m}}(H^d_{\mathfrak{m}}(M)) = 0$ for all $i \not=d$.
	\end{itemize}
\end{corollary}

\begin{proof}
	First prove (i)$ \Longrightarrow $(ii). To this end recall that 
	$d = \grade(\mathfrak{m},M) = \depth_R M$. Then the conclusion follows by 
	\ref{thm-3}. 
	
	For the proof of (ii) $ \Longrightarrow $ (i) we first prove that it will be 
	enough to show the claim for $M = \hat{M}$ the $\mathfrak{m}$-adic completion 
	of $M$ over the completed ring $\hat{R}$. This follows because of 
	$\Lambda_i^{\hat{\mathfrak{m}}}(H^d_{\hat{\mathfrak{m}}}(\hat{M})) \cong 
	\Lambda_i^{\mathfrak{m}}(H^d_{\mathfrak{m}}(M))$ as easily seen (see 
	e.g. \cite[9.8.3]{SS}). Because $R$ is the quotient of a Gorenstein ring 
	the canonical module $K_M$ of $M$ exists and $H^d_{\mathfrak{m}}(M) \cong 
	D(K_M)$ (see \cite{Sp1} resp. \cite[10.3]{SS} for further details). 
	Moreover  
	\[
	\Lambda_i^{\mathfrak{m}}(H^d_{\mathfrak{m}}(M)) \cong
	\Lambda_i^{\mathfrak{m}}(D(K_M)) \cong D(H^i_{\mathfrak{m}}(K_M))
	\] 
	for all $i$ (see e.g. \cite[9.2.5]{SS}). By the assumptions  it follows 
	$H^i_{\mathfrak{m}}(K_M) = 0$ for 
	all $i \not= d$ and $D(H^d_{\mathfrak{m}}(K_M)) =M$. 
	Therefore $K_M$ is a Cohen-Macaulay module as well as $K_{K_M} \cong 
	D(H^d_{\mathfrak{m}}(K_M)) \cong M$, which completes the proof 
	since the canonical module of a Cohen-Macaulay module is Cohen-Macaulay. 
\end{proof}

In \ref{cor-4} (b) it is not enough to assume the vanishing of 
$\Lambda_i^{\mathfrak{m}}(H^d_{\mathfrak{m}}(M))$ for all $i \not= d$. 
This implies (in the case of a complete local ring) that $K_M$ is a 
Cohen-Macaulay module. A finitely generated $R$-module $M$ such that 
$K_{\hat{M}} = \Hom_R(H^d_{\mathfrak{m}}(M),E)$ is a Cohen-Macaulay 
$\hat{R}$-module is called canonically Cohen-Macaulay. See \cite{Sp4} 
for a study of canonically Cohen-Macaulay modules. The following Corollary \ref{cor-5} 
provides an intrinsic characterization of them. The proof is clear by the proof 
of \ref{cor-4}.

\begin{corollary} \label{cor-5}
	With the previous notation the following conditions are equivalent:
	\begin{itemize}
		\item[(i)] $M$ is a canonically Cohen-Macaulay module.
		\item[(ii)] $\Lambda_i^{\mathfrak{m}}(H^d_{\mathfrak{m}}(M)) = 0$ for all $i \not=d$.
	\end{itemize}
\end{corollary}

\section{A Generalization}
At first we shall prove a slight improvement of an aspect of the so-called Greenlees-May 
duality. It was originally proved by Greenlees and May (see \cite{GM}) and by Lipman et al. 
(see \cite{All1}). It was generalized by Porta, Shaul and Yekutieli (see \cite[Theorem 7.14]{PSY} 
and the correction \cite[Theorem 9]{PSY1}).

\begin{proposition} \label{GM}
	Let $R$ denote a commutative ring (not necessarily Notherian). Let $\xx = x_1,\ldots,x_m$ and $\yy = 
	y_1,\ldots,y_n$ denote two weakly proregular system of elements and $I = (\xx)R$ and $J = (\yy)R$. Suppose that 
	$\Rad J \subseteq \Rad I$. Then in the derived category there is an isomorphism 
	\[
	{\rm{R}} \Hom_R(\Rgam (M),{\rm{R}}\Gamma_J(N)) \cong \Llam ({\rm{R}} \Hom_R(M,N))
	\]
	for two $R$-modules $M,N$.
\end{proposition}

\begin{proof}
	Let $F_{\cdot} \qism M$ be a free resolution of $M$. Let $\check{L}_{\xx}$ and $\check{L}_{\yy}$ 
	denote the bounded free resolutions of the \v{C}ech complexes $\check{C}_{\xx}$ and 
	$\check{C}_{\yy}$ resp. (see \cite[6.2.2 and 6.2.3]{SS}). Then the following complex 
	\[
	\Hom_R(\check{L}_{\xx} \otimes_R F_{\cdot},\check{L}_{\yy} \otimes_R N) 
	\cong \Hom_R(F_{\cdot},\Hom_R(\check{L}_{\xx}, \check{L}_{\yy} \otimes_R N))
	\]
	is a representative of the left complex in the statement. Now the natural morphism 
	\[
	\Hom_R(\check{L}_{\xx}, \check{L}_{\yy} \otimes_R N) \to \Hom_R(\check{L}_{\xx},N)
	\]
	is a quasi-isomorphism since $\Rad \yy R \subseteq \Rad \xx R$ (see \cite[6.5.4]{SS}). 
	So there is a quasi-isomorphism 
	\[
	\Hom_R(F_{\cdot},\Hom_R(\check{L}_{\xx}, \check{L}_{\yy} \otimes_R N)) \to \Hom_R(F_{\cdot},\Hom_R(\check{L}_{\xx},N)) \cong \Hom_R(\check{L}_{\xx},\Hom_R(F_{\cdot},N)).
	\]
	Finally recall that the last complex in the previous sequence is a representative of the second 
	complex in the statement (see again \cite{SS}).
\end{proof}

As an application of \ref{GM} there is an alternative proof of \ref{cor-2} and an extension to the relative 
situation. 

\begin{corollary}\label{cor-6}
	Let $I, J$ denote two ideals in a Noetherian ring $R$ with 	$\Rad J \subseteq \Rad I$. 
	Let $M,N$ be two finitely generated $R$-modules. Suppose that $H_I^i(M) = 0$ for 
	all $i \not= g$ and $H^j_J(N) = 0$ for all $j \not= h$. Then there are 
	isomorphisms 
	\[
	\Ext_R^{i+g-h}(H^g_I(M),H^h_J(N)) \cong \hat{R}^I \otimes_R\Ext_R^i(M,N)
	\]
	for all $i \geq 0$.
\end{corollary}

\begin{proof}
	We apply the isomorphism of \ref{GM}. Under the assumptions the $i$-th cohomology 
	of the left complex is $\Ext_R^{i+g-h}(H^g_I(M),H^h_J(N))$. Let $F_{\cdot} \qism M$ 
	be a free resolution by finitely generated free $R$-modules. 
	The complex at the right in \ref{GM} is represented by 
	\[
	\Hom_R(\check{L}_{\xx},\Hom_R(F_{\cdot},N)) \cong \Hom_R(F_{\cdot},\Hom_R(\check{L}_{\xx},N)) \cong \Hom_R(F_{\cdot},R) \otimes_R \Hom_R(\check{L}_{\xx},N)
	\]
	(see \cite[11.1.2]{SS}). Since $N$ is finitely generated $\Hom_R(\check{L}_{\xx},N)$ 
	is quasi-isomorphic to $\hat{R}^I \otimes_RN$. Therefore there is a quasi-isomorphism 
	of the previous complexes  to 
	\[
	\Hom_R(F_{\cdot},R)\otimes_R N \otimes_R \hat{R}^I \cong \Hom_R(F_{\cdot},N) \otimes_R 
	\hat{R}^I.
	\]
	Now note that $\hat{R}^I$ is $R$-flat and taking cohomology provides the claim. 
\end{proof}

For the isomorphism $\Llam ({\rm{R}} \Hom_R(M,N)) \cong \hat{R}^I \otimes_R {\rm{R}} \Hom_R(M,N)$ for finitely generated $R$-modules $M,N$ see also \cite[2.7]{Fa} resp. \cite{SS}. Moreover, note that $\grade(\Ann_RM,N) = \inf \{ i \in \mathbb{N}| \Ext_R(M,N) \not= 0\}$ (see e.g. \cite[1.2.10]{BrH}). 

\begin{corollary} \label{cor-7}
	Suppose that the ideals $I,J \subset R$ satisfy the assumptions of \ref{cor-6}. 
	Let $M,N$ denote two finitely generated $R$-modules with $c = \cd (I,M)$ and 
	$h = \grade (J,N)$. Then there is an isomorphism 
	$\Hom_R(H^c_I(M),H^h_J(N)) \cong \hat{R}^I \otimes_R \Ext_R^{h-c}(M,N)$. 
\end{corollary}

\begin{proof} 
	We have $c = \sup\{i \in \mathbb{N}| H_I^i(M) \not= 0\}$ and $h = \inf \{H^i_J(N) 
	\not= 0\}$. By view of 
	\cite[2.1]{Fhb} it follows that $H^i({\rm{R}} \Hom_R(\Rgam (M),{\rm{R}}\Gamma_J(N))) = 
	0$ for all $i < -c+h$ and  
	\[
	H^{-c+h}({\rm{R}} \Hom_R(\Rgam (M),{\rm{R}}\Gamma_J(N))) \cong 
	\Hom_R(H^c_I(M),H^h_J(N)).
	\]
	This finishes the proof as in \ref{cor-6}.
\end{proof}

\section{Examples}

\begin{example} \label{expl-1}
	(A)  With the notation of \ref{thm-3} there is an exact sequence 
	\[
	0 \to H_1(\Hom_R(\check{L}_{\xx},C_M^{\cdot}(I))) \to 
	\Lambda_g^I(H^g_I(M)) \stackrel{\tau}{\longrightarrow} \hat{M}^I \to 
	H_0(\Hom_R(\check{L}_{\xx},C_M^{\cdot}(I))) \to \Lambda_{g-1}^I(H^g_I(M))\to 0.
	\]
	It would be of some interest to have an intrinsic description of 
	$\Ker \tau$ and $\Coker \tau$.
	Here $\xx = x_1,\ldots,x_t$ denotes a system of elements that generates 
	the maximal ideal $\mathfrak{m}$ up to the radical, e.g. a system of parameters of $R$. \\
	(B) Let $M$ denote a $d$-dimensional module with $d > g = \depth_R M >1$. 
	Suppose 
	that $H^i_{\mathfrak{m}}(M)$ is finitely generated for all $i \not= d$, 
	i.e. $M$ is a generalized Cohen-Macaulay module. Then 
	$\Lambda_i^{\mathfrak{m}}(H^g_{\mathfrak{m}}(M)) = 0$ for all $i \not =0$. 	
	For the 
	homology of 
	$\Hom_R(\check{L}_{\xx},C_M^{\cdot}(I))$ there is the spectral sequence
	\[
	E^2_{i,j} = \Lambda_i^{\mathfrak{m}}(H^j_{\mathfrak{m}}(M)) 
	\Longrightarrow E^{\infty}_{i-j} = \Lambda_{i-j}^{\mathfrak{m}}(C^{\cdot}_M(\mathfrak{m})).
	\]
	Since $H^i_{\mathfrak{m}}(M)$ is finitely generated for $i \not=d$ 
	it follows that $\Lambda_i^{\mathfrak{m}}(H^j_{\mathfrak{m}}(M)) = 0$ for $j 
	\not= d$ and $i > 0$. Whence there is a  partial degeneration to 
	$\Lambda_0^{\mathfrak{m}}(C^{\cdot}_M(\mathfrak{m})) \cong 
	\Lambda_d^{\mathfrak{m}}(H^d_{\mathfrak{m}}(M))$ at the stage $i-j = 0$ 
	because $\depth_R M > 1$. Therefore $M \cong 
	\Lambda_d^{\mathfrak{m}}(H^d_{\mathfrak{m}}(M))$. Because $R$ is the quotient 
	of a Gorenstein ring 
	the canonical module $K_M$ of $M$ exists and $H^d_{\mathfrak{m}}(M) \cong 
	D(K_M)$ (see \cite{Sp1} resp. \cite[10.3]{SS} for further details). 
	Moreover $\Lambda_d^{\mathfrak{m}}(D(K_M)) \cong D(H^d_{\mathfrak{m}}(K_M)) 
	\cong K_{K_M} \cong M$. But this is well-known under the assumptions about $M$ 
	(see \cite{Sp1}). 
\end{example}

\begin{example} \label{expl-2}
	Let $R = \Bbbk[[x_1,x_2,x_3]]$ the formal power series ring in three variables 
	over the field $\Bbbk$. Let $I = (x_1,x_2,x_2x_3,x_1x_3) 
	= (x_1,x_2)\cap (x_2,x_3) \cap (x_1,x_3) \supset J = (x_1x_2,x_1x_3+x_2x_3)$.
	Then $\Rad I = \Rad J$ as easily seen. Because $J$ is a complete intersection 
	of $\height (J) = 2$ it follows $H^i_J(R) =0$ for $i \not= 2$ and therefore 
	$R \cong 	\Hom_R(H^g_J(R), H^g_J(R)) \cong \Hom_R(H^g_I(R), H^g_I(R)).$
\end{example}

\begin{example} \label{expl-3}
	As a consequence of  Hartshorne's Second Vanishing Theorem (see \cite{Hr2}) 
	one has the following: 
	Let $C \subset \mathbb{P}^3_{\Bbbk}$ denote a connected curve. Then 
	$H^i_{I}(R) = 0$ for $i \not= 2$, where $I \subset R = 
	\Bbbk[[x_0,x_1,x_2,x_3]]$ is the saturated defining ideal of $C$. That is 
	$R \cong \Hom_R(H^2_I(R),H^2_I(R))$.
	
	Note that the image $\mathbb{P}^1_{\Bbbk} \to \mathbb{P}^3_{\Bbbk}, 
	(s:t) \mapsto (s^4:s^3t:st^3:t^4)$, is set-theoretically a complete 
	intersection for $\Char(\Bbbk) = p >0$ (see \cite{Hr3} or \cite{RS}) 
	while this is open in the case of $\Char(\Bbbk) = 0$. So, it seems that 
	the endomorphism ring $\Hom_R(H^2_I(R),H^2_I(R))$ is not sensitive 
	enough in order to solve this problem.
\end{example}

Let $I \subset R$ denote an ideal of a local ring $(R,\mathfrak{m})$ with 
$\grade(I,R) = g$. In general the structure of the endomorphism ring $\Hom_R(H^g_I(R),H^i_I(R))$ 
is difficult to understand. In general it is not a finitely generated $R$-module. 
The following example was suggested by R. Hartshorne (see \cite{Hr4}).

\begin{example} \label{expl-4} (see \cite[{\S} 3]{Hr4} and \cite[Example 3.6]{Sp6})
	Let $\Bbbk$ denote a field and $R = \Bbbk[[x,y,u,v]]/(xv-yu),$ where
	$\Bbbk[[x,y,u,v]]$ denotes the formal power series ring in four variables over
	$\Bbbk$. Let $I = (x,y)R.$ Then $R$ is a Gorenstein ring with $\dim R = 3, 
	\dim R/I = 2$ and $\grade (I,R) = 1.$
	It follows that $H^i_I(R) = 0$ for $i \not= 1,2.$ Moreover $\Supp H^2_I(R) 
	\subset \{\mathfrak m\}.$ We use the truncation
	complex  (see \ref{def-1}).
	By the long exact sequence of the local cohomology it induces a short exact sequence 
	\[
	0 \to H^2_I(R) \to H^2_{\mathfrak m}(H^1_I(R)) \to E \to 0.
	\]
	(see \cite[Lemma 2.2]{Sp7} for the details).
	R. Hartshorne (cf. \cite[{\S} 3]{Hr4}) has shown that the socle of
	$H^2_I(R)$ is not a finite dimensional $\Bbbk$-vector space. Therefore,
	the socle of $H^2_{\mathfrak m}(H^1_I(R))$ is infinite dimensional. Moreover 
	there are the following isomorphisms
	\[
	\Hom_R(H^1_I(R), H^1_I(R)) \cong \Ext^1_R(H^1_I(R), R) \cong \Hom_R(H^2_{\mathfrak m}(H^1_I(R)), E),
	\]
	where the second isomorphism follows by the Local Duality Theorem 
	(see e.g. \cite[10.4.2]{SS} for not necessarily finitely generated modules).
	By the Nakayama Lemma this means that $\Hom_R(H^1_I(R), H^1_I(R))$ is not a 
	finitely generated $R$-module.
	
	A bit more is true: Let $S = \Bbbk[[u,v,a]]$ and $R \to S$ the 
	ring homomorphism induced by $x \mapsto au, y \mapsto av$ which is injective. 
	In \cite[Example 3.6]{Sp6} it is shown that $S \cong \Hom_R(H^1_I(R), 
	H^1_I(R))$. This proves that $ \Hom_R(H^1_I(R), H^1_I(R))$ is a Noetherian 
	ring but not finitely generated over $R$. Also it yields another proof that 
	the socle dimension of $H^2_I(R)$ is not finite since $S$ is not a finitely 
	generated $R$-module.
\end{example}

{\bf Acknowledgement.} The author is deeply grateful to the reviewer for essential 
comments and suggestions for further generalizations of the author's approach. 

\bibliographystyle{siam}

\bibliography{wiegand}

\end{document}